\documentclass[12pt]{amsart}
\usepackage{amssymb, amstext, amscd, amsmath}
\usepackage{amsthm}
\usepackage{txfonts}
\usepackage{graphicx}
\newtheorem{thm}{Theorem}[section]

\newtheorem{lem}[thm]{Lemma}

\theoremstyle{definition}

\numberwithin{equation}{section}

\begin{document}

\title[A proof of Lickorish and Wallace's theorem]
{A proof of Lickorish and Wallace's theorem}

\author{Qiang E}
\address{School of Mathematical Sciences, Dalian University of
Technology, Dalian 116024, China} \email{eqiangdut@gmail.com}

\author{Fengchun Lei}
\address{School of Mathematical Sciences, Dalian University of
Technology, Dalian 116024, China} \email{ffcclei@yahoo.com.cn}

\author{Fengling Li}
\address{School of Mathematical Sciences, Dalian University of
Technology, Dalian 116024, China} \email{fenglingli@yahoo.com.cn}
\thanks{This research is supported by NSFC No.11101058 }



\begin{abstract}In this paper, we give a simple proof of Lickorish and Wallace's theorem, which states that every closed orientable 3-manifold is obtained by surgery on some link in $S^3$.
\end{abstract}
\maketitle

\section{Introduction}
Let $K$ be a knot embedded in an orientable closed 3-manifold $M$. Remove $\eta(K)$ which
denotes an open regular neighborhood of $K$, and glue a solid torus
$V$ back to obtain a new 3-manifold $M'$, such that the meridian
curve of $\partial V$ is identified with \emph{L} which is an essential,
unoriented simple closed curve on $\partial(M\setminus\eta(K))$. $M'$ is said to be obtained by \emph{Dehn surgery on $K$ along $L$ in $M$}.
Let $K^*\subset M'$ be a core loop of $V$. We call $K^*$ the dual knot of $K$ in $M'$.
We remark that if some Dehn surgery on $K$ in $M$ produces a 3-manifold $M'$,
then $K^*$ admits a Dehn surgery yielding $M$.

Let $\emph{M}$ be a connected orientable closed 3-manifold. If there
is a closed surface $\emph{S}$ which cuts $\emph{M}$ into two
handlebodies $\emph{V}$ and $\emph{W}$ with $S=\partial V=\partial
W$, then we say $\emph{M}$ has a \emph{Heegaard splitting}, denoted
by $M=V\cup_{S}W$; and $\emph{S}$ is called a \emph{Heegaard
surface} of $\emph{M}$. The \emph{Heegaard genus} of a 3-manifold $M$ is defined
as the genus of the minimal genus Heegaard surface and denoted by $g(M)$. It is
well known that every compact orientable 3-manifold admits a Heegaard splitting.

The following result is known as the fundamental theorem of surgery the-
ory. It has been proven independent by Lickorish and Wallace. Lickorish's
proof is based on the Lickorish twist theorem, which states that any ori-
entation preserving homeomorphism is isotopic to a composition of a finite
number of Dehn twists. A different approach to Lickorish's proof, due to A.
Wallace, may be found in\cite{AH}. Its starting point is the fact that any closed,
orientable 3-manifold W is the boundary of a smooth, compact, orientable
4-manifold.

\begin{thm}\cite{WB}\cite{AH}
Every closed orientable 3-manifold $M$ is obtained by surgery on some link in $S^{3}$.
\end{thm}

In this paper, we give a simple proof of Theorem 1.1. First we show
the following lemma, which is the main result of this paper.

\begin{lem}
Let $M$ be a closed orientable 3-manifold which is not homomorphic
to $S^3$. Then there exists a link in $M$ such that surgery on this
link in $M$ may produce a 3-manifold $M'$ such that $g(M')<g(M)$.
\end{lem}

Theorem 1.1 follows from it immediately.
\begin{proof}(of Theorem 1.1 under the Lemma 1.2)
We prove the theorem by induction on $g(M)$. It is trivial when
$g(M)=0$, since $M\cong S^3$. Suppose the theorem holds when
$g(M)< g$. Let $M$ be a closed 3-manifold which admit a genus $g$
minimal Heegaard splitting. By Lemma 1.2, There exits a link $L$
in $M$ such that surgery on this link may produce a 3-manifld $M'$
such that $g(M')<g$. By our assumption, There exits a link $L'$ in
$M'$ such that surgery on this link may produce $S^{3}$. We may
suppose  $L$ and $L'$ is disjoint in $M$. Thus surgery on the
union link in $M$ may produce $S^3$. $M$ could be obtained by
surgery on the dual link in $S^3$.

\end{proof}

\section{Preliminaries}
Terms which are not defined in this section are all standard, refer to \cite{MS} \cite{DR} and \cite{WJ}.

Let $M=V\cup_{S}W$ be a Heegaard splitting. It is said to be \emph{stabilized} if there are two
essential disks $D_{1}\subset V$ and $D_{2}\subset W$ such that
$|\partial D_{1}\cap\partial D_{2}|=1$; otherwise, it is
\emph{unstabilized}. A genus $g$ Heegaard splitting $V\cup_{S}W$ is
stabilized if and only if there exists a genus $(g-1)$ Heegaard
splitting $V'\cup_{S'}W'$ such that $V\cup_{S}W$ is obtained from
$V'\cup_{S'}W'$ by adding an unknotted handle, and then we say
$V\cup_{S}W$ is obtained from $V'\cup_{S'}W'$ by a stabilization.

Suppose $g(S)>1$. The \emph{(Hempel) distance}\cite{JH} between two essential simple closed curves $\alpha$ and $\beta$ in $S$, denoted
by $d(\alpha, \beta)$, is the smallest integer $n \geq 0$ such that there
is a sequence of essential simple closed curves $\alpha=\alpha_{0},
\alpha_{1}...,\alpha_{n}=\beta$ in $S$ such that $\alpha_{i-1}$ is
disjoint from $\alpha_{i}$ for $1\leq i\leq n$. The \emph{(Hempel) distance} of the
Heegaard splitting $V\cup_{S}W$ is $d(S)=Min \{ d(\alpha, \beta)
\},$ where $\alpha$ bounds an essential disk $D$ in $V$ and $\beta$ bounds an essential disk $E$ in $W$. Moreover, if $M$ is closed, each $\alpha_i$ could be chosen so that it is non-separating in $S$, see\cite{FM}.

A simple closed curve $c$ embedded in the boundary surface of a handlebody $H$ is called \emph{primitive} in $H$,
if it intersects an essential disk of $H$  in a single point. If we push $c$ into the interior of $H$ then the copy is a core loop of $H$, say $c'$. Notice that doing any surgery on $c'$ yields a handlebody. Now consider a spanning annulus $A$ in the compression body $H\setminus\eta(c')$ with
$c\subset\partial A$, and the other boundary component of $A$, says $c''$, lies in $\partial N(c')$. Surgery on $c'$ along $c''$ in $H$ will yield a handlebody $H'$ where $c$ bounds a meridian disk in $H'$. For the sake of convenience, in the following statement, we denote that $H'$ is obtained by a \emph{P-M surgery on $c$ in $H$} (See Figure 1).
\begin{figure}[htbp]
    \includegraphics[width=7.5cm]{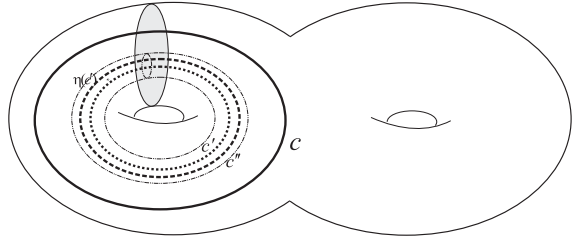}\\
  \caption{P-M surgery on $c$: Dehn surgery on $c'$ along $c''$}\label{fg}
\end{figure}
\section{A proof of Lemma 1.2}

\begin{proof}
Let $V\cup_F W$ be a Heegaard splitting of $M$ such that
$g(F)=g(M)$.

\textbf{Claim.} There exists a finite sequence of essential simple closed curves in $F$, say $[\alpha_{1},...,\alpha_{m}]$ , such that $\alpha_{1}$ is primitive in $V$, $\alpha_{m}$ is primitive in $W$, and  $\alpha_{i-1}$ intersects $\alpha_{i}$ transversely in a single point for $2\leq i\leq m$.

\emph{Proof.} It is easy to check when $g(F)=1$, using the Farey
graph. So we assume that $g(F)>1$ and $d(F)=n$, then there is a
sequence of non-separating simple closed curves $a_0, a_1, ... ,
a_n$ in $F$, such that $a_0$ bounds a disk in $V$, $a_n$ bounds a
disk in $W$ and
 $a_i\cap a_{i-1}=\emptyset$, for $1\leq i\leq n$. For each pair $a_{i}$ and $a_{i+1}$, we may find a simple closed curve $b_{i}$ such that $|a_{i}\cap b_{i}|=|b_{i}\cap a_{i+1}|=1$ for $0\leq i\leq n-1$ (See Figure 2). Let $\alpha_{2i+1}=b_i$ for $0\leq i\leq n-1$ and $\alpha_{2j}=a_j$ for $1\leq j\leq n-1$, and let $m=2n-1$. Then $\{\alpha_k\}$ for $1\leq k\leq m$ satisfies our claim.   $\quad\square$
 \begin{figure}[htbp]
    \includegraphics[width=7.5cm]{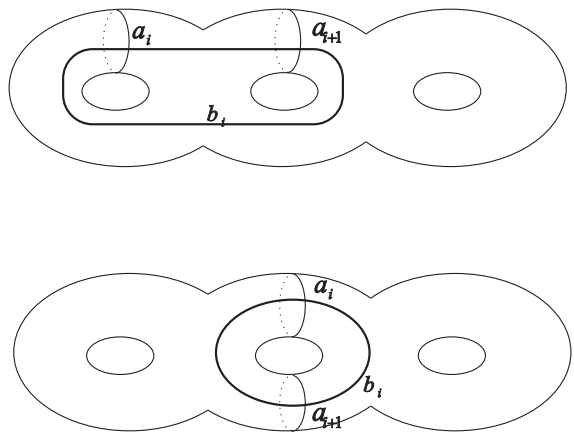}\\
  \caption{$b_i$ which intersects both $a_i$ and $_{i+1}$ in a single point}\label{bo}
\end{figure}

Let $V$ be $V_0$ and $V_i$ be the handlebody obtained from P-M
surgery on $\alpha_i$ in $V_{i-1}$.  Such P-M surgery could last
because $\alpha_i$ is primitive in $V_{i-1}$ and is a meridian of
$V_i$ for $1\leq i\leq m$. Repeating this process allows us to
produce $M'=V_m\cup_{F}W$. If we consider $V$ as $(F\times I)\cup
V'$, where $V'$ is a copy of $V$ and $\partial V=F=F\times \{1\}$,
we may push $\alpha_i$ so that $\alpha_i'$ lies in $F\times
\{i/n\}$, where $n$ is a sufficiently large number. Then $M'$ is
obtained by a surgery on $\alpha_{1}'\cup \alpha_{2}'...\cup
\alpha_{m }'$ in $M$. Since $\alpha_m$ is a meridian of $V_m$ and
primitive in $W$, $V_m\cup W$ is a stabilized Heegaard splitting
with genus $g(M)$. So we have $g(M')<g(M)$.
\end{proof}


\bibliographystyle{amsplain}

\end{document}